\documentclass[12pt]{amsart}
\usepackage{amsmath, amssymb, euscript, hyperref}

\newtheorem{thm}{Theorem}[section]
\newtheorem{lem}[thm]{Lemma}
\newtheorem{prop}[thm]{Proposition}
\theoremstyle{definition}

\newtheorem{rem}[thm]{Remark}
\newtheorem{cor}[thm]{Corollary}

\begin{document}

\bibliographystyle{amsplain}

\address{Azer Akhmedov, Department of Mathematics,
North Dakota State University,
Fargo, ND, 58102, USA}
\email{azer.akhmedov@ndsu.edu}

 \address{Michael P. Cohen, Department of Mathematics,
North Dakota State University,
Fargo, ND, 58102, USA}
\email{michael.cohen@ndsu.edu}

\bigskip

\begin{center}{\bf Some applications of H\"older's theorem in groups of analytic diffeomorphisms of 1-manifolds} \end{center}

\bigskip

\begin{center} Azer Akhmedov and Michael P. Cohen \end{center}

\subjclass[2010]{}

\bigskip

ABSTRACT: {\small We obtain a simple obstruction to embedding groups into the analytic diffeomorphism groups of 1-manifolds. Using this, we classify all RAAGs which embed into $\mathrm{Diff}_{+}^{\omega }(\mathbb{S}^1)$. We also prove that a branch group does not embed into $\mathrm{Diff}_{+}^{\omega }(\mathbb{S}^1)$.}

\section{Introduction}

 \medskip

 It is a classical fact (essentially due to H\"older, [cf] \cite{N1}) that if a subgroup $\Gamma \leq \mathrm{Homeo}_{+}(\mathbb{R})$ acts freely then $\Gamma $ is Abelian. This is obtained by first showing that $\Gamma $ is necessarily left-orderable and Archimedean, and then showing that Archimedean groups are Abelian. 
 
 \medskip
 
 It is interesting to ask the same question for an arbitrary manifold, i.e. if $M$ is an orientable manifold and a group $\Gamma \leq \mathrm{Homeo}_{+}(M)$ acts freely then is $\Gamma $ necessarily Abelian? Motivated by this question, we will say that an orientable manifold $M$ is {\em H\"older} if any freely acting subgroup of $\mathrm{Homeo}_{+}(M)$ is Abelian. Similarly, if $M$ is an orientable manifold with boundary $\partial M$, then we will say that $M$ is {\em H\"older} if any freely acting subgroup $\Gamma \leq \mathrm{Homeo}_{+}(M, \partial M)$ is Abelian. At the moment we have very little understanding of H\"older manifolds. For example, we do not even know if the 2-dimensional cube $I^2$ is H\"older.\footnote{This question seems to be related to the question of Calegari and Rolfsen in \cite{CR} about left-orderability of $\mathrm{Homeo}_{+}(I^n, \partial I^n), n\geq 2$. Also, in \cite{C}, Calegari proves a striking topological consequence of the freeness but only for $\mathbb{Z}^2$-actions, namely, for a free $C^1$ class $\mathbb{Z}^2$ action by the orientation preserving homeomorphisms of the plane, the Euler class of the action must vanish!}  On the other hand, it follows from the Lefshetz fixed-point theorem that even-dimensional spheres $\mathbb{S}^{2n}, n\geq 1$ (and more generally, even-dimensional rational homology spheres) are trivially  H\"older. Interestingly, the circle $\mathbb{S}^1$ is also a (non-trivial) H\"older manifold, that is any freely acting subgroup of orientation preserving circle homeomorphisms is necessarily Abelian (see \cite{N1}). 

 \medskip
 
 H\"older's theorem provides strong and sometimes surprising obstructions to embedding groups into the homeomorphism group of the circle. For example, it immediately implies that (as remarked in \cite{N1}) a finitely generated infinite torsion group does not embed in $\mathrm{Homeo}_{+}(\mathbb{S}^1)$. Let us emphasize that even for the simplest 2-manifolds such as $\mathbb{S}^2$ or $\mathbb{T}^2$ this question is still unsettled; recently, some interesting partial results have been obtained by N.Guelman and I.Liousse (see \cite{GL1} and \cite {GL2}). 
 
 \medskip
  
 In this paper we study obstructions for embeddability of groups into the group of analytic diffeomorphisms of 1-manifolds. Throughout the paper, $M^1$ will denote a compact oriented 1-manifold, so we will assume that either $M^1 = I$ or $M^1 = \mathbb{S}^1$, with a fixed orientation. $\mathrm{Diff}_{+}^{\omega }(M^1)$ will denote the group of orientation preserving analytic diffeomorphisms of $M^1$.

 \medskip

 The following lemma shows that the commutativity relation among the elements of $\mathrm{Diff}_{+}^{\omega }(I)$ is transitive.

 \medskip

 \begin{lem}\label{lem:infectious} Let $f, g, h\in \mathrm{Diff}_{+}^{\omega }(I)$ be non-trivial elements such that $f$ commutes with $g$ and $g$ commutes with $h$. Then $f$ commutes with $h$.
 \end{lem}

 \medskip

 {\bf Proof.} A non-trivial analytic diffeomorphism has finitely many fixed points. Therefore if two non-trivial analytic diffeomorphisms $\phi _1, \phi _2$ commute then $Fix(\phi _1) = Fix(\phi _2)$.

 \medskip

  Thus we obtain that $Fix(f) = Fix(g) = Fix(h)$. Let us assume that $Fix(f) = Fix(g) = Fix(h) = \{p_1, \dots , p_n\}$, where, $0 = p_1 < p_2 < \dots < p_n = 1$.

  \medskip

 Now, assume that $f$ and $h$ do not commute, and let $H$ be a subgroup of $\mathrm{Diff}_{+}^{\omega }(I)$ generated by $f$ and $h$. $H$ fixes all the points $p_1, \ldots , p_n$. Then, by H\"older's Theorem, $H$ contains a non-trivial element $\phi $ such that $\phi $ has a fixed point $p\in (0,1)$ distinct from $p_1, \ldots , p_n$. Let $p$ lie in between $p_{i}$ and $p_{i+1}$ for some $i\in \{1,\dots ,n\}$. Since $[\phi , g] = 1$, we obtain that $\phi $ must have infinitely many fixed points in $(p_{i}, p_{i+1})$. Contradiction. $\square $

 \bigskip

  Lemma \ref{lem:infectious} implies strong restrictions for the embedding of groups into Diff$_{+}^{\omega }(M^1)$. We would like to observe the following fact, which will be useful for our consideration of right-angled Artin groups in Section 2.

  \medskip

  \begin{lem}\label{lem:K_N} Let $f_1, \dots , f_N$ be non-identity elements of $\mathrm{Diff}_{+}^{\omega }(I)$. Let also $G$ be a simple graph on $N$ vertices $v_1, \dots , v_N$ such that $(v_i v_j)$ is an edge whenever $f_i$ and $f_j$ commute. Assume that $G$ is connected. Then $G$ is isomorphic to a complete graph $K_N$.
 \end{lem}

  \medskip

  {\bf Proof.} Let $1\leq i < j\leq N$. We want to show that $v_i$ and $v_j$ are connected with an edge. Let $(u_0, u_1, \dots , u_m)$ be a path in $G$ such that $u_0 = v_i, u_m = v_j$. If $(u_0, u_j)$ is an edge, and $j < m$, then, by Lemma \ref{lem:infectious}, $(u_0, u_{j+1})$ is an edge. Since $(u_0,u_1)$ is an edge, by induction, we obtain that $(u_0,u_m)$ is an edge. $\square $

 \medskip

 Since H\"older's Theorem holds for the analytic diffeomorphism group of the circle as well, slightly weaker versions of both of the lemmas \ref{lem:infectious} and \ref{lem:K_N} generalize to $\mathrm{Diff}_{+}^{\omega }(M^1)$.

 \medskip
 
 \begin{lem} \label{lem:circle} Let $f, g, h\in \mathrm{Diff}_{+}^{\omega }(\mathbb{S}^1)$ be elements of infinite order, $f$ commutes with $g$ and $g$ commutes with $h$. Then $f^m$ commutes with $h^m$ for some $m\geq 1$.
 \end{lem}
   
 \medskip
 
 {\bf Proof.} Let $\Gamma $ be a subgroup generated by $f$ and $h$. If $\Gamma $ acts freely then by H\"older's Theorem it is Abelian, and the claim follows. If $\Gamma $ does not act freely then let $p$ be a fixed point for some non-trivial $\gamma \in \Gamma $. By analyticity, $\gamma$ has finitely many fixed points, and therefore since $[g,\gamma]=1$, $p$ is a fixed point of $g^n$ for some $n\geq 1$. Then $g^n$ has finitely many fixed points, and $[f,g^n]=[h,g^n]=1$, so $f^{i}$ and $h^{j}$ both fix $p$ for some positive integers $i,j$. Setting $m=ij$, $p$ is a fixed point of both $f^m$ and $h^m$. The claim now follows from Lemma \ref{lem:infectious} applied to $f^m$, $g^n$, and $h^m$, since the subgroup of elements of $\mathrm{Diff}_+^\omega(\mathbb{S}^1)$ which fix $p$ may be identified with a subgroup of $\mathrm{Diff}_+^\omega(I)$. $\square$
 

 
 Now we obtain an analogue of Lemma \ref{lem:K_N}. 
 
 \medskip
 
 \begin{lem}\label{lem:analogue}   Let $f_1, \dots , f_N$ be elements of $\mathrm{Diff}_{+}^{\omega }(\mathbb{S}^1)$ of infinite order. Let also $G$ be a simple graph on $N$ vertices $v_1, \dots , v_N$ such that $(v_i v_j)$ is an edge whenever $f_i^m$ and $f_j^m$ commute for some $m\geq 1$. Assume that $G$ is connected. Then $G$ is isomorphic to a complete graph $K_N$. $\square $
\end{lem}
    
  \begin{rem} It is easy to see that Lemma \ref{lem:infectious} and Lemma \ref{lem:circle} both fail for the the group of orientation preserving analytic diffeomorphisms of the real line. Indeed, for a real number $a\neq 0$, let $f_a(x) = \frac{1}{2}\sin(ax) + x$ and $g_a(x) = x+\frac{2\pi }{a}$. Then, for all distinct $a, b\in \mathbb{R}\backslash \{0\}$, we have $[f_a, g_a] = 1$ and $[g_a, g_b] = 1$ while for some $a, b$, the diffeomorphisms $f_a^n, f_b^n$ do not commute for any $n\geq 1$.  However, it is possible to prove a similar but much weaker statement for $\mathrm{Diff}_{+}^{\omega }(\mathbb{R})$, and then an analogue of Lemma \ref{lem:analogue}, at the expense of increased technicality. We skip the relevant discussion here.    
  \end{rem}
    
\bigskip

 \section{Classification of RAAGs in $\mathrm{Diff}_{+}^{\omega }(M^1)$}

 \medskip

  A right angled Artin group (RAAG) is defined as follows. Let $G = (V, E)$ be a finite simple graph with a set of vertices $V = \{v_1, \ldots , v_n\}$, and $\Gamma $ be a finitely presented group given by the presentation $$\Gamma = \langle x_1, \ldots , x_n \ | \ [x_i,x_j] = 1 \ \mathrm{iff} \  (v_iv_j) \ \mathrm{is \ an \ edge \ in} \  G \rangle .$$ Thus the graph $G$ defines the group $\Gamma $. RAAGs have been studied extensively in the past decades from combinatorial, algebraic and geometric points of view. In the recent paper \cite{BKK} the authors prove the following very interesting result.

  \begin{thm} Every RAAG embeds in $\mathrm{Diff}_{+}^{\infty }(M^1)$.
  \end{thm}

 The following theorem gives a necessary and sufficient condition under which a RAAG embeds in $\mathrm{Diff}_{+}^{\omega }(M^1)$.

 \begin{thm} A RAAG $\Gamma = \Gamma (G)$ embeds in $\mathrm{Diff}_{+}^{\omega }(M^1)$ if and only if every connected component of the graph $G$ is a complete graph.
 \end{thm}

 \medskip

  {\bf Proof.} We will present a constructive proof. The ``only if" part follows immediately from Lemmas \ref{lem:K_N} and \ref{lem:circle} for $M^1=I$ and $M^1=\mathbb{S}^1$, respectively. For the ``if" part, let $C_1, \dots , C_m$ be all connected components of $G$, such that $C_i$ is a complete graph on $m_i$ vertices.

  \medskip

  We let $f, g_n, n\geq 1$ be non-trivial orientation preserving analytic diffeomorphisms of $M^1$ satisfying the following conditions:

 \medskip

  (c1) for every finite subset $A\subseteq \mathbb{N}$ the diffeomorphisms $g_n, n\in A$ generate a free Abelian group of rank $|A|$.

 \medskip

  (c2) the subgroup generated $f$ and $g_n, n\geq 1$ is isomorphic to the free product $\langle f \rangle \ast \langle g_1, g_2, \dots \rangle $.

 \medskip

  Now, let $h_{i,n} = f^ig_nf^{-i}$ for all non-negative $i, n$. For all $i\in \{1,\dots ,m\}$, let $H_i$ be a subgroup generated by $h_{i, 1}, \dots , h_{i, m_i}$.

  \medskip

  Notice that, for all $i\in \{1, \dots , m\}$,  the group $H_i$ is a free Abelian group of rank $m_i$. Let $\Gamma $ be a subgroup generated by $H_1, \dots , H_m$. Then $\Gamma = H_1\ast H_2\ast \dots \ast H_m$ thus $\Gamma $ is isomorphic to the RAAG of the graph $G$.

  \medskip

  It remains to show that there exist elements $f, g_n, n\geq 1$ satisfying conditions (c1) and (c2). Let $\alpha _1, \alpha _2, \dots $ be rationally independent numbers in $(1, \pi )$ such that $\log \alpha _1, \log \alpha _2 , \dots $ are also rationally independent. In the case of $M^1 = I$, we choose $g_n(x) = \frac{\alpha _n x}{(\alpha _n-1)x+1}$ for all $n\geq 1$, and in the case of $M^1 = \mathbb{S}^1$, we choose $g_n$ to be the orientation preserving rotation by the angle $\alpha_n$.
  
  \medskip

   Then $g_n, n\geq 0$ satisfy condition (c1). Now, we need to prove that there exists $f\in \mathrm{Diff}_{+}^{\omega }(M^1)$ such that  condition (c2) holds.

  \medskip

  We indicate the construction for $M^1=I$.  Let $F$ be a free group formally generated by letters $f, g_1, g_2, \dots $ (so $F$ is a free group with infinite rank; also, we abuse the notation by denoting the generators of $F$ by $g_1, g_2, \dots $ which are already defined analytic diffeomorphisms, and by $f$ which is an analytic diffeomorphism we intend to define). Let $A = \{f, f^{-1}, g_1, g_1^{-1}, g_2, g_2^{-1}, \dots \}$, and $W_1, W_2, \ldots $ be {\em all} non-trivial reduced words in the alphabet $A$ such that each $W_i$ contains $f$ or $f^{-1}$. 

  \medskip

  Let $D$ be an open connected and bounded domain in $\mathbb{C}$ such that $D$ contains the real interval $[0,1]$. We build $f$ as a sum $f=\displaystyle\sum_{n=0}^\infty \omega _n$ of analytic functions $\omega _n$ on $D$.  Intuitively speaking, each $n$-th summand $\omega_n$ induces a small perturbation which prevents the corresponding $n$-th word $W_n(f,g_1,g_2,...)$ from reducing to identity.  We recursively define the maps $\omega_0,\omega_1,\omega_,\dots$, as well as sequences $D_1,D_2,\dots$ and $\epsilon_1,\epsilon_2,\dots$ of real numbers as follows.
  
 \medskip

 Let $\omega _0$ be an identity map, i.e. $\omega _0(z) = z, \forall z\in D$. Fix an arbitrary point $p\in(0,1)$. We let $\omega _1:D\to \mathbb{C}$ be an analytic function such that

 (i) $\omega_1(0)=\omega_1(1)=0$,
 
 (ii) $\omega_1(x)$ is real for all $x\in\mathbb{R}\cap D$,
 
 (iii) $\displaystyle\sup_{z\in D}|\omega_1(z)|<\frac{1}{4}$,
 
 (iv) $\displaystyle\sup_{x\in[0,1]}|\omega_1'(x)|<\frac{1}{4}$, and
 
 (v) $W_1(f_1,g_1,g_2,\dots)(p) \neq p$, where $f_1=\omega_0+\omega_1$.
 
 \medskip
 
 Let $D_1=|W_1(f_1,g_1,g_2,\dots)(p)-p|$.  Let $\epsilon_1<1$ be a positive number so small that if $\displaystyle\sup_{x\in[0,1]}|\varphi(x)-f_1(x)|<\epsilon_1$, then $$\displaystyle\sup_{x\in[0,1]}|W_1(\varphi,g_1,g_2,\dots)(x)-W_1(f_1,g_1,g_2,\dots)(x)|<\frac{D_1}{2}.$$
 
 Suppose now $m\geq 2$ and $\omega_1, \dots, \omega_{m-1}$ are chosen so that the quantity $D_{m-1} = |W_{m-1}(f_{m-1},g_1,g_2,\dots)(p)-p|$ is strictly positive where $f_{m-1}=\omega_0+\omega_1+\dots+\omega_{m-1}$.  Also suppose $\epsilon_1,\dots,\epsilon_{m-1}$ are chosen so that $\epsilon_{m-1}\leq\dots\leq\epsilon_1\leq 1$, and if $\displaystyle\sup_{x\in[0,1]}|\varphi(x)-f_{m-1}(x)|<\epsilon_{m-1}$, then $$\displaystyle\sup_{x\in[0,1]}|W_i(\varphi,g_1,g_2,\dots)(x)-W_i(f_{m-1},g_1,g_2,\dots)(x)|<\frac{D_i}{2}$$ for each $1\leq i\leq m-1$.  Then choose $\omega_m:D\to\mathbb{C}$ such that:
 
 \medskip
 
 (i) $\omega_m(0)=\omega_m(1)=0$,

 (ii) $\omega_m(x)$ is real for all $x\in\mathbb{R}\cap D$,
 
 (iii) $\displaystyle\sup_{z\in D}|\omega_m(z)|<\frac{\epsilon_{m-1}}{2^{m+1}}$,
 
 (iv) $\displaystyle\sup_{x\in[0,1]}|\omega_m'(x)|<\frac{1}{2^{m+1}}$, and
 
 (v) $W_m(f_m,g_1,g_2,\dots)(p) \neq p$, where $f_m=\omega_0+\omega_1+\dots+\omega_m$.
 

 \medskip

 Let $f(z)=\displaystyle\sum_{n=0}^\infty\omega_n(z)$ for $z\in D$.  By condition (iii) $f$ is a uniform limit of analytic functions on an open bounded domain and hence analytic.  Moreover, $f'(x)$ is real for all $x\in[0,1]$, and $f'(x)=\omega_0'(x)+\displaystyle\sum_{n=1}^\infty\omega_n'(x)\geq 1-\displaystyle\sum_{n=1}^\infty\frac{1}{2^{n+1}}>0$ by condition (iv), so $f$ is increasing and hence an analytic diffeomorphism of $I$ since it fixes $0$ and $1$.  For every word $W_n$, we have $W_n(f,g_1,g_2,\dots)(p)\neq p$, since $\displaystyle\sup_{x\in[0,1]}|f(x)-f_n(x)|\leq\displaystyle\sum_{i={n+1}}^\infty\frac{\epsilon_i}{2^{i+2}}\leq\frac{\epsilon_n}{2}<\epsilon_n$ and therefore 
 
 \begin{align*}
 |W_n(f,g_1,g_2,\dots)(p)-p| &\geq |W_n(f_n,g_1,g_2\dots)(p)-p|-|W_n(f,g_1,g_2\dots)(p)\\
 & \hspace{1cm} -W_n(f_n,g_1,g_2\dots)(p)|\\
 &> D_n-\frac{D_n}{2}=\frac{D_n}{2}. 
 \end{align*}
 
 \medskip
 
 Thus condition (c2) is satisfied.  (Obtaining an analytic diffeomorphism of the circle satisfying (c2) may be done similarly.)  $\square$

 \bigskip

 \section{More Applications}

  In this section we will observe other consequences of Lemma \ref{lem:K_N} and Lemma \ref{lem:analogue}.

  \begin{prop} \label{prop:centralizer} Let $\Gamma $ be a group with a sequence of elements $g_1, g_2, \dots $ such that the centralizer $C_{g_i}$ is a proper subset of the centralizer $C_{g_{i+1}}$ for every $i\geq 1$. Then $\Gamma $ does not embed in $\mathrm{Diff}_{+}^{\omega }(I)$.
  \end{prop}

  {\bf Proof.} We will prove a much stronger fact. By the assumption, in $\Gamma $, there exist elements $g_1, g_2, h_1, h_2$ such that $h_1$ belongs to the centralizer of $g_1$, $h_1, h_2$ belong to the centralizer of $g_2$ while $h_2$ does not belong to the centralizer of $g_1$.

  \medskip

  Let $V = \{g_1, g_2, h_1, h_2\}$. Consider a simple graph $G = (V, E)$ with a vertex set $V$ and where the edge set $E$ defined as follows: $(v_iv_j)\in E$ iff $v_i$ and $v_j$ are distinct and $[v_i,v_j] = 1$ in the group $\Gamma $.

  \medskip

  Notice that $(g_1,h_1), (g_2, h_1), (g_2, h_2)$ are edges in $G$ thus $G$ is connected. By Lemma \ref{lem:K_N} $G$ is isomorphic to $K_4$. Hence $[g_1, h_2]=1$. Contradiction. $\square $

 \medskip

 The existence of strictly increasing chain of centralizers is an interesting property for groups; most notably in the area of branch groups. 

 \medskip

 \begin{prop} A branch group does not embed in $\mathrm{Diff}_{+}^{\omega }(M^1)$.
 \end{prop}

\medskip

 {\bf Proof.} a) For $M^1 = I$, this immediate from Lemma \ref {lem:infectious}. Indeed, by definition of a branch group (see \cite {BGS}), it contains elements $f, g, h$ where $[f,g]=[g,h]=1$ while $[f,h]\neq 1$.

 \medskip
 
  b) For $M^1 = \mathbb{S}^1$, let $\Gamma \leq \mathrm{Diff}_{+}^{\omega }(M^1)$ be a branch group. Then, by definition of a branch group, $\Gamma $ has a finite index subgroup of the form $\Gamma _1\times \Gamma _2\times \dots \times \Gamma _n$ where $n\geq 2$ and none of the groups $\Gamma _i, 1\leq i\leq n$ is virtually Abelian. By H\"older's Theorem for the circle, there exists a non-trivial diffeomorphism $f\in \Gamma _n$ such that $Fix(f)\neq \emptyset $. Let $p_1, \dots , p_N$ be all fixed points of $f$ listed along the orientation of the circle $\mathbb{S}^1$.

 \medskip
 
 Let also $H = \{\gamma \in \Gamma _1 \ | \ Fix(\gamma ) \neq \emptyset \}$. Then for every $\gamma \in H\backslash \{1\}$, by analyticity, we have $Fix(\gamma ) = Fix(f)$. On the other hand, if $\gamma (p_i) = p_j$ for some $\gamma \in \Gamma _1$ and $i, j \in \{1, \dots , N\}$, then $\gamma (p_{i+1}) = p_{j+1}$. Hence $H$ is a finite index subgroup of $\Gamma _1$. Also, for every $h_1, h_2\in H$, we have $[h_1, f] = [h_2, f] = 1$. Then, by Lemma \ref{lem:infectious}, $[h_1, h_2] = 1$. Thus $H$ is Abelian. But then $\Gamma _1$ is virtually Abelian. Contradiction. $\square $ 

 \medskip
 
 Now, we would like to observe several other corollaries of the results from Section 1. 
 
 \medskip
 
 \begin{cor} \label{cor:center} A non-Abelian group with a non-trivial center does not embed in $\mathrm{Diff}_{+}^{\omega }(I)$.
 \end{cor}

 \medskip

 This fact (with modified statement) can be deduced already for the group $\mathrm{Diff}_{+}^{2}(I)$ using Kopell's Lemma. More precisely, it follows from Kopell's Lemma that an irreducible subgroup of $\mathrm{Diff}_{+}^{2}(I)$ has a trivial center. 

 \medskip
 
 Notice that $SL(n,\mathbb{Z}), n\geq 3$ contains a copy of the integral $3\times 3$ unipotent subgroup. Then, by Margulis arithmeticity result of higher rank lattices, and by Corollary \ref{cor:center} and Lemma \ref{lem:circle}, we obtain the following corollary. 

 \begin{cor} \label{cor:lattice} A lattice in $SL(n, \mathbb{R}), n\geq 3$ does not embed in $\mathrm{Diff}_{+}^{\omega }(M^1)$.
 \end{cor}

\medskip

 Let us emphasize that it is already known that a finite index subgroup of $SL(n, \mathbb{Z}), n\geq 3$ does not embed in $\mathrm{Homeo}_{+}(I)$, see \cite{W}; for an arbitrary lattice in $SL(n, \mathbb{R}), n\geq 3$ this question is still open. On the other hand, it is also known that a lattice in $SL(n, \mathbb{R}), n\geq 3$ does not embed in $\mathrm{Diff}_{+}(M^1)$ as proved by Ghys and Burger-Monod, and it is also known that any infinite discrete group with property $(T)$ does not embed in $\mathrm{Diff}_{+}^{1+\alpha }(\mathbb{S}^1), \alpha > \frac{1}{2}$ as proved by Navas \cite{N2}. (in particular, the result of Corollary \ref{cor:lattice} is not new.) 
 
 \medskip
 
 It is interesting to study representations of general Artin groups into the group of analytic diffeomorphisms of manifolds. Notice that, by Corollary \ref{cor:center} and Lemma \ref{lem:circle}, the braid groups on $n\geq 3$ strings (another special subclass of Artin groups) do not embed in $\mathrm{Diff}_{+}^{\omega }(M^1)$.

\end{document}